\documentclass[12pt]{amsart}
\usepackage{amsmath}
\usepackage{hyperref}
\usepackage{amssymb}
\usepackage{latexsym}
\usepackage{amscd}
\usepackage{graphicx} %We can use any other package if necessary
\usepackage{amsthm}
\usepackage{mathrsfs}
\usepackage{xypic}
\usepackage{bm}

\newdimen\AAdi%
\newbox\AAbo%
\def\AAk#1#2{\s_etbox\AAbo=\hbox{#2}\AAdi=\wd\AAbo\kern#1\AAdi{}}%
\def\AAr#1#2#3{\s_etbox\AAbo=\hbox{#2}\AAdi=\ht\AAbo\raise#1\AAdi\hbox{#3}}%
\font\tenmsb=msbm10 at 12pt \font\sevenmsb=msbm7 at 8pt
\font\fivemsb=msbm5 at 6pt
\newfam\msbfam
\textfont\msbfam=\tenmsb \scriptfont\msbfam=\sevenmsb
\scriptscriptfont\msbfam=\fivemsb
\def\Bbb#1{{\tenmsb\fam\msbfam#1}}
\newtheorem{thm}{Theorem}[section]
\newtheorem{lem}{Lemma}[section]
\newtheorem{cor}{Corollary}[section]
\newtheorem{rem}{Remark}[section]
\newtheorem{pro}{Proposition}[section]

\newcommand{\ba}{\begin{array}}
\newcommand{\ea}{\end{array}}

\newcommand{\Section}[2]{\setcounter{equation}{0}
\allowdisplaybreaks
\section[#1]{#2}}

\def\n{\nabla}

\def\ir#1{\mathbb R^{#1}}

\def\f#1#2{\frac{#1}{#2}}

\def\grs#1#2{\bold G_{#1,#2}}

\def\a{\alpha}
\def\be{\beta}

\def\p#1{\partial #1}

\def\de{\delta}
\def\De{\Delta}

\def\ep{\varepsilon}
\def\eps{\epsilon}
\def\G{\Gamma}
\def\g{\gamma}

\def\la{\lambda}
\def\La{\Lambda}
\def\om{\omega}
\def\Om{\Omega}
\def\th{\theta}

\def\Si{\Sigma}

\def\w{\wedge}

\def\Hess{\mbox{Hess}}
\def\R{\Bbb{R}}

\def\tr{\mbox{tr}}
\def\U{\Bbb{U}}
\def\lan{\langle}
\def\ran{\rangle}
\def\ra{\rightarrow}

\def\aint#1{-\hskip -4.5mm\int_{#1}}


\begin{document}
\title
[The Gauss image and Bernstein type theorems] {The Gauss image of 
entire graphs of higher codimension and Bernstein type theorems}

\author
[J. Jost, Y. L. Xin and Ling Yang]{J. Jost, Y. L. Xin and Ling
Yang}
\address{Max Planck Institute for Mathematics in the
Sciences, Inselstr. 22, 04103 Leipzig, Germany.}
\email{jost@mis.mpg.de}
\address {Institute of Mathematics, Fudan University,
Shanghai 200433, China.} \email{ylxin@fudan.edu.cn}
\address{Max Planck Institute for Mathematics in the
Sciences, Inselstr. 22, 04103 Leipzig, Germany.}
\email{lingyang@mis.mpg.de}
\thanks{The second named author is grateful to the Max Planck
Institute for Mathematics in the Sciences in Leipzig for its
hospitality and  continuous support. He is also partially
supported by NSFC and SFMEC }

\begin{abstract}
Under suitable conditions on the range of the Gauss map of a 
complete submanifold of Euclidean space with parallel mean 
curvature, we construct  a strongly subharmonic function and  
derive a-priori estimates for the harmonic Gauss map. The 
required conditions here are more general than in previous work 
and they therefore enable us to improve substantially previous 
results for  the Lawson-Osseman problem concerning the regularity 
of minimal submanifolds in higher codimension and to derive 
Bernstein type results.
\end{abstract}
\maketitle

\Section{Introduction}{Introduction}

% This is a continuation of the second part of our previous paper 
% \cite{j-x-y} which deals with Bernstein problem of hypersurfaces 
% in Euclidean space via  a-priori estimates of harmonic maps into 
% the unit sphere.

We consider an oriented $n$-dimensional submanifold $M$ in
$\ir{n+m}$ with $n\ge 3,\; m\ge 2.$  The 
Gauss map $\g:M\to \grs{n}{m}$ maps $M$ into a Grassmann
manifold. In fact,  for codimension $m=1$, this Grassmann
manifold $\grs{n}{1}$ is  the unit sphere $S^n$. In this paper,
however, we are interested in the case $m\ge 2$ where the geometry of
this  Grassmann manifold is more complicated. By the theorem of Ruh-Vilms \cite{r-v}, $\g$ is
harmonic if and only if $M$ has parallel mean curvature. This result
applies thus in particular to the case where $M$ is a {\it minimal}
submanifold of Euclidean space. 
% On the
% other hand, the situation is much more complicated in higher
% codimension, as was revealed by Lawson-Osserman \cite{l-o}.

Now, the Bernstein problem for entire minimal graphs is one of 
the central problems in geometric analysis. Let us summarize the 
status of this problem, first for the case of codimension 1. The 
central result is that an  entire minimal graph $M$ of dimension 
$n\le 7 $ and codimension 1  has to be planar, but there are 
counterexamples to such a Bernstein type theorem in dimension $8$ 
or higher. However, when the additional condition is imposed that 
the slope of the graph be uniformly bounded, then a theorem of  
Moser \cite{m}, called a weak  Bernstein theorem,  asserts that 
such an $M$ in arbitrary dimension has to be planar. Thus, the 
counterexamples arise from a non-uniform behavior at infinity. In 
fact, by a general scaling argument, the Bernstein theorems are 
intimately related to the regularity question for the minimal 
hypersurface equation. 

A natural and important question then is to what extent such 
Bernstein type theorems generalize to entire minimal graphs of 
codimension $m\ge 2$. Moser's result has been extended to  higher 
codimension by Chern-Osserman for dimension $n=2$ \cite{c-o} and 
Barbosa \cite{b} and Fisher-Colbrie \cite{fc} for dimension 
$n=3$. For dimension $n=4$ and codimension $m=3$, however, there 
is a counterexample given by Lawson-Osserman \cite{l-o}. In fact, 
their paper emphasizes the stark contrast between the cases of 
codimension 1 and greater than 1 for the minimal submanifold 
system, concerning regularity, uniqueness, and existence. The 
Lawson-Osserman problem then is concerned with a systematic 
understanding of the analytic aspects of the minimal submanifold 
system in higher codimension. As in the case of codimension 1, 
the Bernstein problem provides a  key towards this aim. 

While the work of Lawson-Osserman produced a counterexample for a
general Bernstein theorem, there are also some positive results 
in this direction which we shall now summarize. 
Hildebrandt-Jost-Widman \cite{h-j-w} started a systematic 
approach on the basis of the aforementioned Ruh-Vilms theorem. 
That is, they developed and  employed the theory of harmonic maps 
and the convex geometry of Grassmann manifolds, and obtained 
Bernstein type results in general dimension and codimension. 
Their main result says that a Bernstein result holds if the image 
of the Gauss map is contained in a strictly convex distance ball. 
Since the Riemannian sectional curvature of $\grs{n}{m}$ is 
nonnegative, the maximal radius of such a convex ball is bounded. 
In codimension 1, this in particular reproduces Moser's theorem, 
and in this sense, their result is optimal. For higher 
codimension, their result can be improved, for the following 
reason. Since the sectional curvature of $\grs{n}{m}$ for $n,m 
\ge 2$ is not constant, there exist larger convex sets than 
geodesic distance balls, and it turns out that harmonic (e.g. 
Gauss) maps with values in such convex sets can still be well 
enough controlled. In this sense, the results of  \cite{h-j-w} 
could be improved  by Jost-Xin \cite{j-x}, Wang \cite{wang} and 
Xin-Yang \cite{x-y1}. In \cite{j-x}, the largest such 
geodesically convex set in a Grassmann manifold
 was found. 

Formulating it somewhat differently, the harmonic map approach is based on the fact that the
 composition of a harmonic map with a convex function is a subharmonic
 function, and by using quantitative estimates for such subharmonic
 functions, regularity and Liouville type results for harmonic maps
 can be obtained. The most natural such convex function is the squared
 distance from some point, when its domain is restricted to a suitably
 small ball. As mentioned, the largest such ball on which a squared
 distance function is convex was utilized in \cite{h-j-w}. As also 
 mentioned, however, this result is not yet optimal, and other    convex functions
were systematically utilized in \cite{x-y1}. In that paper, also the
fundamental connection between estimates for the second fundamental
form of minimal submanifolds and estimates for their Gauss maps was
systematically explored. On this basis, the fundamental curvature
estimate technique, as developed by Schoen-Simon-Yau \cite{s-s-y}
and Ecker-Huisken \cite{e-h}, could be used in \cite{x-y1}.

Still, there remains a large quantitative gap between those positive
results and the counterexample of Lawson-Osserman. In this situation,
it could either be that Bernstein theorems can be found under more
general conditions, or that there exist other counterexamples in the
so far unexplored range. 

In the present paper, we make a step towards closing this gap in 
the positive direction. We identify a geometrically natural 
function $v$ on a Grassmann manifold and a natural quantitative 
condition under which the precomposition of this function with a 
harmonic (Gauss) map is (strongly) subharmonic (Theorem 
\ref{thm1}). When the precomposition of $v$ with the Gauss map of 
a complete minimal submanifold is bounded, then that submanifold 
is an entire graph of bounded slope. On one hand, this is the 
first systematic example in harmonic map regularity theory where 
this auxiliary function is not necessarily convex. On the other 
hand, the Lawson-Osserman's counterexample can also be readily 
characterized in terms of this function. Still, the range of 
values for $v$ where we can apply our scheme is strictly 
separated from the value of $v$ in that example. Therefore, still 
some gap remains which should be explored in future work. 

Our work also finds its natural position in the general regularity
theory for harmonic maps. Also, once we have a strongly 
subharmonic function, we could derive  Bernstein type results 
within the frame work of geometric measure theory, by the 
standard blow-down procedure and appeal to Allard's regularity 
theorem \cite{a}. By building upon the work of many people on 
harmonic map regularity, we can obtain more insight, however. In 
particular, we shall use  the iteration method of \cite{h-j-w}, 
we can explore the relation with curvature estimates,  and we 
shall utilize a version of the telescoping trick (Theorem 
\ref{thm2}) to finally obtain a quantitatively controlled Gauss 
image shrinking process (Theorem \ref{thm3} and Theorem 
\ref{thm4}). In this way, we can understand why the submanifold 
is flat as the Bernstein result asserts.  More precisely, we 
obtain the following Bernstein type result, which substantially 
improves our previous results.

\begin{thm}\label{thm5}
Let $z^\a=f^\a(x^1,\cdots,x^n),\ \a=1,\cdots,m$, be smooth
functions defined everywhere in $\R^n$ ($n\geq 3,m\geq 2$).
Suppose their graph $M=(x,f(x))$ is a submanifold with parallel
mean curvature in $\R^{n+m}$. Suppose that there exists a number
$\be_0<3$ such that
\begin{equation}
\De_f=\Big[\det\Big(\de_{ij}+\sum_\a \f{\p f^\a}{\p x^i}\f{\p
f^\a}{\p x^j}\Big)\Big]^{\f{1}{2}}\leq \be_0.\label{be2}
\end{equation}
Then $f^1,\cdots,f^m$ has to be affine linear, i.e., it  represents an
affine $n$-plane.
\end{thm}

The essential point is to show that $v:=\De_f$ is subharmonic when
$<3$. In fact, when  $v \le  \be_0 <3$, then $\Delta v \ge K_0 |B|^2$
where $K_0$ is a positive constant and $B$ is the second fundamental
form of $M$ in $\R^{n+m}$. This principle is not new. Wang \cite{wang}
has given conditions under which $\log v$ is subharmonic and has
derived Bernstein results from this, as indicated above. He only needs that $v$ be
uniformly bounded by some constant, not necessarily $<3$, but in
addition that there exist some $\delta >0$ such that for any two
eigenvalues $\la_i, \la_j$ with $i\neq j$, the inequality $|\la_i
\la_j|\le 1 -\delta$ holds (the latter condition means in geometric
terms that $df$ is strictly area decreasing on any two-dimensional subspace). Since subharmonicity of $\log v$ is a
weaker property than subharmonicity of $v$ itself, his computation is
substantially easier than ours, and our results cannot be deduced from
his. In fact, $v^2=\prod (1+\la_i^2)$, and while the condition of
\cite{j-x} which can be reformulated as $v^2$ being bounded away from
4 implies the condition of \cite{wang} so that the latter result
generalizes the former, the condition needed in the present paper is
only 
the weaker one that $v^2$ be bounded away from 9. 

In fact,  somewhat more refined results can be obtained, as will 
be pointed out in the final remarks of this paper.

\Section{Geometry of Grassmann manifolds}{Geometry of Grassmann manifolds}\label{s1}

Let $\R^{n+m}$ be an $(n+m)$-dimensional Euclidean space. Its 
oriented $n$-subspaces constitute the Grassmann manifold
$\grs{n}{m}$, which is the Riemannian symmetric space of compact
type $SO(n+m)/SO(n)\times SO(m).$

$\grs{n}{m}$ can be viewed as a submanifold of some  Euclidean space
via the Pl\"ucker embedding. The restriction of the 
Euclidean inner product on $M$ is denoted by $w:\grs{n}{m}\times \grs{n}{m}\ra \R$
$$w(P,Q)=\lan e_1\w\cdots\w e_n,f_1\w\cdots\w f_n\ran=\det W$$
where $P$ is spanned by a unit $n$-vector $e_1\w\cdots\w e_n$, $Q$ is spanned by another unit $n$-vector $f_1\w\cdots
\w f_n$, and $W=\big(\lan e_i,f_j\ran\big)$. It is well-known that
$$W^T W=O^T \La O$$
with $O$ an orthogonal matrix and
$$\La=\left(\begin{array}{ccc}
            \mu_1^2 &   &  \\
                    & \ddots &  \\
                    &        & \mu_n^2
            \end{array}\right).$$
Here each $0\leq \mu_i^2\leq 1$. Putting $p:=\min\{m,n\}$, then 
at most $p$ elements in $\{\mu_1^2,\cdots, \mu_n^2\}$ are not 
equal to $1$. Without loss of generality,  we can assume 
$\mu_i^2=1$ whenever $i>p$. We also note that the $\mu_i^2$ can 
be expressed as
\begin{equation}\label{di1a}
\mu_i^2=\frac{1}{1+\la_i^2}. \end{equation} 

The Jordan angles between $P$ and $Q$ are defined by
$$\th_i=\arccos(\mu_i)\qquad 1\leq i\leq p.$$
The distance between $P$ and $Q$ is defined by
\begin{equation}\label{di}
d(P, Q)=\sqrt{\sum\th_i^2}. \end{equation} 
Thus, (\ref{di1a}) becomes  
\begin{equation}\label{di2}
\la_i=\tan\th_i. \end{equation}

In the sequel, we shall 
assume $n\geq m$ without loss of generality. We use the
summation  convention and agree on the ranges of indices:
$$1\leq  i,j,k,l\leq n,\; 1\leq \a,\be,\g\leq m,\;  a, b,\cdots =1,\cdots, n+m.$$

Now we fix $P_0\in \grs{n}{m}.$ We represent it by $ n $ vectors
$\eps_i,$, which are complemented by $ m $ vectors $ \eps_{n+\a} $,
such that $ \{\eps_i, \eps_{n+\a} \} $ form an orthonormal base of $
\ir{m+n} $.

Denote
$$\Bbb{U}:=\{P\in \grs{n}{m},\; w(P,P_0)>0\}.$$
We can span an arbitrary $P\in \Bbb{U}$ by $n$-vectors $f_i$:
$$f_i=\eps_i+Z_{i\a}\eps_{n+\a}.$$
The canonical metric in $\Bbb{U}$ can be described as
\begin{equation}\label{m1}ds^2 = tr (( I_n + ZZ^T )^{-1} dZ (I_m + Z^TZ)^{-1} dZ^T
),\end{equation} where $ Z = (Z_{i \a}) $ is an $ (n \times m)
$-matrix and $ I_n $ (res. $ I_m $) denotes the $ (n\times n)
$-identity (res. $ m \times m $) matrix. It is shown that
(\ref{m1}) can be derived from (\ref{di}) in \cite{x}.

For any $P\in\Bbb{U}$, the Jordan angles between $P$ and $P_0$ 
are defined by $\{\th_i\}$.  Let $E_{i\a}$ be the matrix with $1$ 
in the intersection of row $i$ and column $\a$ and $0$ otherwise. 
Then, $\sec\th_i\sec\th_\a E_{i\a}$ form an orthonormal basis of 
$T_P\grs{n}{m}$ with respect to (\ref{m1}). Denote its dual frame 
by $\om_{i\a}.$

Our fundamental quantity will be
\begin{equation}
v(\cdot, P_0):=w^{-1}(\cdot, P_0) \text{ on }\Bbb{U}.
\end{equation} 
For
arbitrary $P\in \U$ determined by an $n\times m$ matrix $Z$, it
is easily seen that
\begin{equation}\label{v}
v(P,P_0)=\big[\det(I_n+ZZ^T)\big]^{\f{1}{2}}=\prod_{\a=1}^m 
\sec\th_\a =\prod_{\a=1}^m \frac{1}{\mu_\alpha}.
\end{equation}
where $\th_1,\cdots,\th_m$ denote the Jordan angles between $P$
and $P_0$.

In this terminology, Hess$(v(\cdot, P_0)$ has been estimated in
\cite{x-y1}. By (3.8) in \cite{x-y1}, we have
\begin{eqnarray}\label{He}\aligned
\Hess(v(\cdot,P_0))&=\sum_{i\neq \a}v\ \om_{i\a}^2+\sum_\a
(1+2\la_\a^2)v\ \om_{\a\a}^2
+\sum_{\a\neq\be} \la_\a\la_\be v(\om_{\a\a}\otimes \om_{\be\be}+\om_{\a\be}\otimes\om_{\be\a})\\
&=\sum_{m+1\leq i\leq n,\a}v\
\om_{i\a}^2+\sum_{\a}(1+2\la_\a^2)v\ \om_{\a\a}^2
                                          +\sum_{\a\neq \be}\la_\a\la_\be v\ \om_{\a\a}\otimes\om_{\be\be}\\
&\qquad\qquad+\sum_{\a<\be}\Big[(1+\la_\a\la_\be)v\Big(\f{\sqrt{2}}{2}(\om_{\a\be}
+\om_{\be\a})\Big)^2\\
&\hskip2in+(1-\la_\a\la_\be)v\Big(\f{\sqrt{2}}{2}(\om_{\a\be}-\om_{\be\a})\Big)^2\Big].
\endaligned
\end{eqnarray}
It follows that
\begin{equation}\label{hess}
 v(\cdot,P_0)^{-1}\Hess(v(\cdot,P_0))
 =g+\sum_\a 2\la_\a^2 \om_{\a\a}^2+\sum_{\a\neq \be}\la_\a\la_\be(\om_{\a\a}\otimes \om_{\be\be}+
\om_{\a\be}\otimes \om_{\be\a}).
\end{equation}

The canonical Riemannian metric on $\grs{n}{m}$ derived from
(\ref{di}) can also be described by the moving frame method. This will
be useful for understanding some of the sequel. Let
$\{e_i,e_{n+\a}\}$ be a local orthonormal frame field in
$\ir{n+m}.$  Let $\{\om_i,\om_{n+\a}\}$ be its dual frame field
so that the Euclidean metric is
$$g=\sum_{i}\om_i^2+\sum_{\a}\om_{n+\a}^2.$$
The Levi-Civita  connection forms $\om_{ab}$ of $\ir{n+m}$ are
uniquely determined by the equations
$$\aligned
&d\om_{a}=\om_{ab}\wedge\om_b,\cr &\om_{ab}+\om_{ba}=0.
\endaligned$$
It is shown in \cite{x} that the canonical Riemannian metric on
$\grs{n}{m}$ can be written as
\begin{equation}\label{m2}
ds^2=\sum_{i,\ \a}\om_{i\, n+\a}^2. \end{equation}

\Section{Subharmonic  functions}{Subharmonic  functions}

Let $M^m\ra \R^{n+m}$ be an isometric immersion with  second
fundamental form $B.$ Around any point $p\in M$, we  choose an
orthonormal frame field $e_i,\cdots, e_{n+m}$ in $\R^{n+m},$ such
that $\{e_i\}$ are tangent to $M$ and $\{e_{n+\a}\}$ normal to
$M.$  The metric on $M$ is $g=\sum_i \om_i^2.$ We have the structure
equations
\begin{equation}\label{str}
\om_{i\  n+\a}=h_{\a ij}\om_j,
\end{equation}
 where $h_{\a ij}$
are the coefficients of second fundamental form $B$ of $M$ in
$\R^{n+m}.$

Let $0$ be the origin of $\R^{n+m}$, $SO(m+n)$ be the Lie group
consisting of all orthonormal frames $(0;e_i,e_{n+\a})$,
$TF=\big\{(p;e_1,\cdots,e_n):p\in M,e_i\in T_p M,\lan
e_i,e_j\ran=\de_{ij}\big\}$ be the principle bundle of orthonormal
tangent frames over $M$, and
$NF=\big\{(p;e_{n+1},\cdots,e_{n+m}):p\in M,e_{n+\a}\in N_p M\big\}$
be the principle bundle of orthonormal normal frames over $M$. Then
$\bar{\pi}: TF\oplus NF\ra M$ is the projection with fiber
$SO(n)\times SO(m)$.

The Gauss map $\g: M\ra \grs{n}{m}$ is defined by
$$\g(p)=T_p M\in \grs{n}{m}$$
via the parallel translation in $\R^{n+m}$ for every $p\in M$. Then the following  diagram commutes
$$\CD
 TF \oplus NF @>i>> SO (n+m)  \\
 @V{\bar\pi}VV     @VV{\pi}V \\
 M  @>{\g}>>  \grs{n}{m}
\endCD$$
where $i$ denotes the inclusion map and $\pi: SO(n+m)\ra \grs{n}{m}$ is defined by
$$(0;e_i,e_{n+\a})\mapsto e_1\w\cdots\w e_n.$$

It follows that
\begin{equation}\label{edg}
|d\g|^2=\sum_{\a,i,j}h_{\a ij}^2=|B|^2.
\end{equation}

(\ref{hess}) was computed for the metric (\ref{m1}) whose
corresponding coframe field is $\om_{i \a}.$ Since (\ref{m1}) and
(\ref{m2}) are equivalent to each other, at any fixed point
$P\in\grs{n}{m}$ there exists an isotropic group action, i.e., an 
$SO(n)\times SO(m)$ action, such that $\om_{i\a}$ is 
transformed  to $\om_{i\ n+\a}$, namely, there are a local tangent
frame field and a local normal frame field such that at the 
point under consideration, 
\begin{equation}\label{str2}
\om_{i\  n+\a}=\g^*\om_{i \a}.
\end{equation}
In conjunction with (\ref{str}) and (\ref{str2}) we obtain
\begin{equation}\label{hij}
\g^*\om_{i\a}=h_{\a ij}\om_j.
\end{equation}

By the Ruh-Vilms theorem \cite{r-v},  the mean curvature of $M$
is parallel if and only if its Gauss map is a harmonic map. Now,
we assume that $M$ has parallel mean curvature.

We define
\begin{equation}
v:=v(\cdot,P_0)\circ \g,
\end{equation}
This function $ v$   on $M$ will be the source of the basic 
inequality for this paper. Its geometric significance is seen 
from the following observation. If the $v-$ function has an upper 
bound (or the $w-$function has a positive lower bound), $M$ can 
be described as an entire graph on $\ir{n}$ by $f:\ir{n}\to 
\ir{m}$, provided $M$ is complete. In this situation, $\la_i$ is 
the singular values of $df$ and
\begin{equation}\label{v1}
v=\Big[\det\Big(\de_{ij}+\sum_\a \f{\p f^\a}{\p x^i}\f{\p
f^\a}{\p x^j}\Big)\Big]^{\f{1}{2}}
\end{equation}
Using the composition formula, in conjunction with (\ref{hess}),
(\ref{edg}) and (\ref{hij}), and the fact that  $\tau(\g)=0$ (the
tension field of the Gauss map vanishes \cite{r-v}), we deduce the
important 
formula  of Lemma 1.1
in \cite{fc} and Prop. 2.1 in \cite{wang}. 
\begin{pro}Let $M$ be an $n-$submanifold in $\ir{n+m}$ with parallel mean
curvature. Then
\begin{equation}\label{Dv}
\De v=v|B|^2+v\sum_{\a,j}2\la_\a^2h_{\a,\a j}^2 +v\sum_{\a\neq
\be,j}\la_\a\la_\be(h_{\a,\a j}h_{\be,\be j}+h_{\a,\be
j}h_{\be,\a j}),
\end{equation}
where $h_{\a,ij}$ are the coefficients of the second fundamental form
of $M$ in $\ir{n+m}$ (see (\ref{str}).
\end{pro}

A crucial step in this paper is to find a condition which
guarantees the strong subharmonicity of the $v-$ function on $M$.  
More precisely, under a condition on $v$, we shall bound its Laplacian  from below by a 
positive constant times  squared norm of the second 
fundamental form.

Looking at  the expression (\ref{Dv}),  we group its terms according to the
different types of the indices of the coefficients of the second
fundamental form as follows.
\begin{equation}
v^{-1}\De v=
\sum_\a\sum_{i,j>m}h_{\a,ij}^2+\sum_{j>m}I_j+\sum_{j>m,\a<\be}II_{j\a\be}
+\sum_{\a<\be<\g}III_{\a\be\g}+\sum_\a IV_\a
\end{equation} where
\begin{equation}
I_j=\sum_\a(2+2\la_\a^2)h_{\a,\a j}^2+\sum_{\a\neq
\be}\la_\a\la_\be h_{\a,\a j}h_{\be,\be j},
\end{equation}
\begin{equation}
II_{j\a\be}=2h_{\a,\be j}^2+2h_{\be,\a j}^2+2\la_\a\la_\be h_{\a,\be j}h_{\be,\a j},
\end{equation}
\begin{equation}\aligned
III_{\a\be\g}=&2h_{\a,\be\g}^2+2h_{\be,\g\a}^2+2h_{\g,\a\be}^2\\
&+2\la_\a\la_\be h_{\a,\be\g}h_{\be,\g\a}+2\la_\be\la_\g h_{\be,\g\a}h_{\g,\a\be}+2\la_\g\la_\a h_{\g,\a\be}h_{\a,\be\g}
               \endaligned
\end{equation}
and
\begin{equation}\aligned
IV_\a=&(1+2\la_\a^2)h_{\a,\a\a}^2+\sum_{\be\neq \a}\big(h_{\a,\be\be}^2+(2+2\la_\be^2)h_{\be,\be\a}^2\big)\\
      &+\sum_{\be\neq \g}\la_\be\la_\g h_{\be,\be \a}h_{\g,\g \a}+2\sum_{\be\neq \a}\la_\a\la_\be h_{\a,\be\be}h_{\be,\be\a}.
      \endaligned
\end{equation}

It is easily seen that
\begin{equation}\label{es1}
I_j=(\sum_\a \la_\a h_{\a,\a j})^2+\sum_\a (2+\la_\a^2)h_{\a,\a j}^2\geq 2\sum_\a h_{\a,\a j}^2.
\end{equation}

Obviously
\begin{equation}
II_{j\a\be}=\la_\a\la_\be(h_{\a,\be j}+h_{\be,\a j})^2+(2-\la_\a\la_\be)(h_{\a,\be j}^2+h_{\be,\a j}^2).
\end{equation}
$v=\Big(\prod_\a (1+\la_\a^2)\Big)^{\f{1}{2}}$ implies $(1+\la_\a^2)(1+\la_\be^2)\leq v^2$. Assume
$(1+\la_\a^2)(1+\la_\be^2)\equiv C\leq v^2$, then differentiating both sides implies
$$\f{\la_\a d\la_\a}{1+\la_\a^2}+\f{\la_\be d\la_\be}{1+\la_\be^2}=0.$$
Therefore
\begin{equation}
\aligned
d(\la_\a \la_\be)&=\la_\be d\la_\a+\la_\a d\la_\be\\
                 &=\big[\la_\be^2(1+\la_\a^2)-\la_\a^2(1+\la_\be^2)\big]\f{d\la_\a}{\la_\be(1+\la_\a^2)}\\
                 &=(\la_\be^2-\la_\a^2)\f{d\la_\a}{\la_\be(1+\la_\a^2)}.
                 \endaligned
\end{equation}
It follows that $(\la_\a,\la_\be)\mapsto \la_\a\la_\be$ attains its maximum at the point satisfying $\la_\a=\la_\be$,
which is hence $((C^{\f{1}{2}}-1)^{\f{1}{2}},(C^{\f{1}{2}}-1)^{\f{1}{2}})$.
Thus $\la_\a\la_\be\leq C^{\f{1}{2}}-1\leq v-1$ and moreover
\begin{equation}\label{es2}
II_{j\a\be}\geq (3-v)(h_{\a,\be j}^2+h_{\be,\a j}^2).
\end{equation}

\bigskip

\begin{lem}\label{l1}
 $III_{\a\be\g}\geq
(3-v)(h_{\a,\be\g}^2+h_{\be,\g\a}^2+h_{\g,\a\be}^2)$.
\end{lem}

\begin{proof}
It is easily seen that
$$\aligned
III_{\a\be\g}-&(3-v)(h_{\a,\be\g}^2+h_{\be,\g\a}^2+h_{\g,\a\be}^2)\\
=&(\la_\a h_{\a,\be\g}+\la_\be h_{\be,\g\a}+\la_\g
h_{\g,\a\be})^2+(v-1-\la_\a^2)h_{\a,\be\g}^2\\
&\qquad+(v-1-\la_\be^2)h_{\be,\g\a}^2
+(v-1-\la_\g^2)h_{\g,\a\be}^2.\endaligned$$

If $\la_\a^2,\la_\be^2,\la_\g^2\leq v-1$, then
$III_{\a\be\g}-(3-v)(h_{\a,\be\g}^2+h_{\be,\g\a}^2+h_{\g,\a\be}^2)$
is obviously nonnegative definite. Otherwise, we can assume
$\la_\g^2>v-1$ without loss of generality, then
 $(1+\la_\a^2)(1+\la_\be^2)(1+\la_\g^2)\leq v^2$ implies $\la_\a^2<v-1, \la_\be^2<v-1$.

Denote $s=\la_\a h_{\a,\be\g}+\la_\be h_{\be,\g\a}$, then by the Cauchy-Schwarz inequality,
$$\aligned
s^2&=(\la_\a h_{\a,\be\g}+\la_\be h_{\be,\g\a})^2\\
   &=\Big(\f{\la_\a}{\sqrt{v-1-\la_\a^2}}\sqrt{v-1-\la_\a^2}h_{\a,\be\g}+\f{\la_\be}{\sqrt{v-1-\la_\be^2}}
     \sqrt{v-1-\la_\be^2}h_{\be,\g\a}\Big)^2\\
   &\leq\Big(\f{\la_\a^2}{v-1-\la_\a^2}+\f{\la_\be^2}{v-1-\la_\be^2}\Big)\big((v-1-\la_\a^2)h_{\a,\be\g}^2+(v-1-\la_\be^2)h_{\be,\g\a}^2\big)
\endaligned$$
i.e.
\begin{equation}
 (v-1-\la_\a^2)h_{\a,\be\g}^2+(v-1-\la_\be^2)h_{\be,\g\a}^2\geq \Big(\f{\la_\a^2}{v-1-\la_\a^2}+\f{\la_\be^2}{v-1-\la_\be^2}\Big)^{-1}s^2.
\end{equation}
Hence
\begin{equation}\label{ineq3}\aligned
 &III_{\a\be\g}-(3-v)(h_{\a,\be\g}^2+h_{\be,\g\a}^2+h_{\g,\a\be}^2)\\
 \geq& (s+\la_\g h_{\g,\a\be})^2+\Big(\f{\la_\a^2}{v-1-\la_\a^2}+\f{\la_\be^2}{v-1-\la_\be^2}\Big)^{-1}s^2+(v-1-\la_\g^2)h_{\g,\a\be}^2\\
=&\Big[1+\Big(\f{\la_\a^2}{v-1-\la_\a^2}+\f{\la_\be^2}{v-1-\la_\be^2}\Big)^{-1}\Big]s^2+(v-1)h_{\g,\a\be}^2+2\la_\g
sh_{\g,\a\be}.
\endaligned
\end{equation}
It is well known that $ax^2+2bxy+cy^2$ is nonnegative definite if
and only if $a,c\geq 0$ and $ac-b^2\geq 0$. Hence the right hand
side of (\ref{ineq3})
 is nonnegative definite if and only if
\begin{equation}\label{cond}
 (v-1)\Big[1+\Big(\f{\la_\a^2}{v-1-\la_\a^2}+\f{\la_\be^2}{v-1-\la_\be^2}\Big)^{-1}\Big]-\la_\g^2\geq 0
\end{equation}
i.e.
\begin{equation}\label{cond2}
 \f{1}{v-1-\la_\a^2}+\f{1}{v-1-\la_\be^2}+\f{1}{v-1-\la_\g^2}\leq \f{2}{v-1}.
\end{equation}

Denote $x=1+\la_\a^2$, $y=1+\la_\be^2$, $z=1+\la_\g^2$. Let $C$ be a
constant $\leq v^2$, denote
$$\Om=\big\{(x,y,z)\in \R^3:1\leq x,y<v,\; z>v,\; xyz=C\big\}$$
and $f:\Om\ra \R$
$$(x,y,z)\mapsto \f{1}{v-x}+\f{1}{v-y}+\f{1}{v-z}.$$
We claim $f\leq \f{2}{v-1}$ on $\Om$. Then (\ref{cond2}) follows
and hence
$$III_{\a\be\g}-(3-v)(h_{\a,\be\g}^2+h_{\be,\g\a}^2
+h_{\g,\a\be}^2)$$ is nonnegative definite.

We now verify the claim. For arbitrary $\ep>0$, denote
$$f_\ep=\f{1}{v+\ep-x}+\f{1}{v+\ep-y}+\f{1}{v+\ep-z},$$
then $f_\ep$ is obviously a smooth function on
$$\Om_\ep=\big\{(x,y,z)\in \R^3:1\leq x,y\leq v,\; z\geq v+2\ep,\; xyz=C\big\}.$$
The compactness of $\Om_\ep$ implies the existence of
$(x_0,y_0,z_0)\in \Om_\ep$ satisfying
\begin{equation}\label{sup}
 f_\ep(x_0,y_0,z_0)=\sup_{\Om_\ep} f_\ep.
\end{equation}

Fix $x_0$, then (\ref{sup}) implies that for arbitrary $(y,z)\in \R^2$
satisfying $1\leq y\leq v,\;  z\geq v+2\ep$ and $yz=\f{C}{x_0}$,
we have
$$f_{\ep,x_0}(y, z)=\f{1}{v+\ep-y}+\f{1}{v+\ep-z}\leq
\f{1}{v+\ep-y_0}+\f{1}{v+\ep-z_0}.$$ Differentiating both sides
of $yz=\f{C}{x_0}$ yields $\f{dy}{y}+\f{dz}{z}=0.$ Hence
\begin{equation}\aligned
&d\Big(\f{1}{v+\ep-y}+\f{1}{v+\ep-z}\Big)=\f{dy}{(v+\ep-y)^2}+\f{dz}{(v+\ep-z)^2}\\
=&\Big[\f{y}{(v+\ep-y)^2}-\f{z}{(v+\ep-z)^2}\Big]\f{dy}{y}=\f{((v+\ep)^2-yz)(y-z)}{(v+\ep-y)^2(v+\ep-z)^2}\f{dy}{y}.
\endaligned
\end{equation}
It implies that $f_{\ep,x_0}\left(y, \f{C}{yx_0}\right)$ is
decreasing in $y$ and $y_0=1.$ Similarly,  one can derive
$x_0=1$. Therefore
$$\sup_{\Om_\ep} f_\ep=f_\ep(1,1,C)=\f{2}{v+\ep-1}+\f{1}{v+\ep-C}<\f{2}{v+\ep-1}.$$
Note that $f_\ep\ra f$ and $\Om\subset \lim_{\ep\ra
0^+}\Om_\ep$. Hence by letting $\ep\ra 0$ one can obtain $f\leq
\f{2}{v-1}.$

\end{proof}

\bigskip

\begin{lem}\label{l2}
There exists a positive constant $\ep_0$, such that if $v\leq 3$, then
$$IV_\a\geq \ep_0\big(h_{\a,\a\a}^2+\sum_{\be\neq \a}(h_{\a,\be\be}^2+2h_{\be,\be\a}^2)\big).$$
\end{lem}

\begin{proof}

For arbitrary $\ep_0\in [0,1)$, denote $C=1-\ep_0$, then
\begin{equation}\label{ineq4}
\aligned
 &IV_\a-\ep_0\big(h_{\a,\a\a}^2+\sum_{\be\neq \a}(h_{\a,\be\be}^2+2h_{\be,\be\a}^2)\big)\\
=&(\sum_\be \la_\be h_{\be,\be \a})^2+(C+\la_\a^2)h_{\a,\a\a}^2+\sum_{\be\neq \a}\big[Ch_{\a,\be\be}^2+(2C+\la_\be^2)h_{\be,\be\a}^2+2\la_\a\la_\be h_{\a,\be\be}h_{\be,\be\a}\big].
\endaligned
\end{equation}
Obviously
$$\aligned
C\, h_{\a,\be\be}^2&+C^{-1}\la_\a^2\la_\be^2h_{\be,\be\a}^2+2\la_\a\la_\be h_{\a,\be\be}h_{\be,\be\a}\\
&\geq (C^{\f{1}{2}}h_{\a,\be\be}+C^{-\f{1}{2}}\la_\a\la_\be
h_{\be,\be\a})^2\geq 0,\endaligned$$ hence, the third term of the
right hand side of (\ref{ineq4}) satisfies
\begin{equation}\label{ineq2}
Ch_{\a,\be\be}^2+(2C+\la_\be^2)h_{\be,\be\a}^2+2\la_\a\la_\be h_{\a,\be\be}h_{\be,\be\a}\geq (2C+\la_\be^2-C^{-1}\la_\a^2\la_\be^2)h_{\be,\be\a}^2
\end{equation}

If there exist 2 distinct indices $\be,\g\neq \a$ satisfying
$$2C+\la_\be^2-C^{-1}\la_\a^2\la_\be^2\leq 0$$ and
$$2C+\la_\g^2-C^{-1}\la_\a^2\la_\g^2\leq 0,$$ then $\la_\a^2>C$ and
$$\la_\be^2\geq \f{2C^2}{\la_\a^2-C},\qquad \la_\g^2\geq \f{2C^2}{\la_\a^2-C}.$$
It implies
$$(1+\la_\a^2)(1+\la_\be^2)(1+\la_\g^2)\geq \f{(\la_\a^2+1)(\la_\a^2+2C^2-C)^2}{(\la_\a^2-C)^2}.$$
Define $f:x\in (C,+\infty)\mapsto 
\f{(x+1)(x+2C^2-C)^2}{(x-C)^2}$, then a direct calculation shows
$$(\log f)'=\f{1}{x+1}+\f{2}{x+2C^2-C}-\f{2}{x-C}=\f{(x-C(2C+3))(x+C)}{(x+1)(x+2C^2-C)(x-C)}.$$ It follows that
$f(x)\geq f(C(2C+3))=\f{(2C+1)^3}{C+1}$, i.e.
\begin{equation}\label{ineq7}
v^2\geq (1+\la_\a^2)(1+\la_\be^2)(1+\la_\g^2)\geq \f{(2C+1)^3}{C+1}.
\end{equation}
If $C=1$, then $\f{(2C+1)^3}{C+1}=\f{27}{2}>9$; hence there is $\ep_1>0$, once
$\ep_0\leq \ep_1$, then $C=1-\ep_0$ satisfies $\f{(2C+1)^3}{C+1}>9$,
which causes a contradiction to $v^2\leq 9$.

Hence,  one can find an index $\g\neq \a$, such that
\begin{equation}
2C+\la_\be^2-C^{-1}\la_\a^2\la_\be^2> 0\qquad \text{for arbitrary }\be\neq \a,\g.
\end{equation}
Denote $s=\sum_{\be\neq \g}\la_\be h_{\be,\be\a}$, then by using the 
Cauchy-Schwarz inequality,
\begin{equation}\label{ineq5}\aligned
 (C+\la_\a^2)h_{\a,\a\a}^2&+\sum_{\be\neq
 \a,\g}(2C+\la_\be^2-C^{-1}\la_\a^2\la_\be^2)h_{\be,\be\a}^2\\
& \geq \Big(\f{\la_\a^2}{C+\la_\a^2}+\sum_{\be\neq
\a,\g}\f{\la_\be^2}{2C+\la_\be^2-C^{-1}\la_\a^2\la_\be^2}\Big)^{-1}s^2.
\endaligned\end{equation} Substituting (\ref{ineq5}) and (\ref{ineq2})
into (\ref{ineq4}) yields
\begin{equation}\label{ineq6}
\aligned
IV_\a&-\ep_0\big(h_{\a,\a\a}^2+\sum_{\be\neq \a}(h_{\a,\be\be}^2+2h_{\be,\be\a}^2)\big)\\
&\geq (s+\la_\g h_{\g,\g\a})^2+
\Big(\f{\la_\a^2}{C+\la_\a^2}+\sum_{\be\neq \a,\g}\f{\la_\be^2}{2C+\la_\be^2-C^{-1}\la_\a^2\la_\be^2}\Big)^{-1}s^2\\
&\qquad +(2C+\la_\g^2-C^{-1}\la_\a^2\la_\g^2)h_{\g,\g\a}^2\\
&\geq \Big[1+\Big(\f{\la_\a^2}{C+\la_\a^2}+\sum_{\be\neq
\a,\g}\f{\la_\be^2}{2C+\la_\be^2-C^{-1}\la_\a^2\la_\be^2}\Big)^{-1}\Big]s^2\\
&\qquad
+(2C+2\la_\g^2-C^{-1}\la_\a^2\la_\g^2)h_{\g,\g\a}^2+2\la_\g s
h_{\g,\g\a}.
\endaligned
\end{equation}

Note that when $m=2$, $s=\la_\a h_{\a,\a\a}$ and $\sum_{\be\neq
\a,\g}\f{\la_\be^2}{2C+\la_\be^2-C^{-1}\la_\a^2\la_\be^2}=0$.

The right hand side of (\ref{ineq6}) is nonnegative definite if
and only if
\begin{equation}\label{con1}
 2C+2\la_\g^2-C^{-1}\la_\a^2\la_\g^2\geq 0
\end{equation}
and
\begin{equation}\label{con2}
\Big[1+\Big(\f{\la_\a^2}{C+\la_\a^2}+\sum_{\be\neq \a,\g}\f{\la_\be^2}{2C+\la_\be^2-C^{-1}\la_\a^2\la_\be^2}\Big)^{-1}\Big](2C+2\la_\g^2-C^{-1}\la_\a^2\la_\g^2)-\la_\g^2\geq 0.
\end{equation}

Assume $2C+2\la_\g^2-C^{-1}\la_\a^2\la_\g^2< 0$, then
$\la_\a^2>2C$ and $\la_\g^2> \f{2C^2}{\la_\a^2-2C}$, which implies
$(1+\la_\a^2)(1+\la_\g^2)\geq
\f{(\la_\a^2+1)(\la_\a^2+2C(C-1))}{\la_\a^2-2C}$. Define $f:x\in
(2C,+\infty)\mapsto \f{(x+1)(x+2C(C-1))}{x-2C}$, then
$$(\log f)'=\f{1}{x+1}+\f{1}{x+2C(C-1)}-\f{1}{x-2C}=
\f{x^2-4Cx-2C^2(2C-1)}{(x+1)(x+2C(C-1))(x-2C)}.$$ and hence
$$\min f=f\big(C(2+\sqrt{4C+2})\big)=2C^2+2C+1+2C\sqrt{4C+2}.$$
In particular,  when
$C=1,\; \min f= 5+2\sqrt{6}>9$. There exists $\ep_2>0$, such that once
$\ep_0\leq \ep_2$, one can derive $\min f>9$ and moreover $v^2\geq
(1+\la_\a^2)(1+\la_\g^2)>9$, which contradicts $v\leq 3$.
Therefore (\ref{con1}) holds.

If $2C+\la_\g^2-C^{-1}\la_\a^2\la_\g^2\geq 0$, (\ref{con2})
trivially holds.

At last, we consider the situation when there exists $\g, \, 
\g\neq\a$, such that  
$$2C+\la_\g^2-C^{-1}\la_\a^2\la_\g^2<0.$$ In this case,  
(\ref{con2}) is equivalent to
\begin{equation}\label{con}
\f{\la_\a^2}{C+\la_\a^2}+\sum_{\be\neq \a}\f{\la_\be^2}{2C+\la_\be^2-C^{-1}\la_\a^2\la_\be^2}\leq -1.
\end{equation}
Noting that
$$\f{\la_\be^2}{2C+\la_\be^2-C^{-1}\la_\a^2\la_\be^2}=\f{C}{C-\la_\a^2}-\f{2C^3}{(C-\la_\a^2)^2}\f{1}{1+\la_\be^2
+\f{\la_\a^2+C(2C-1)}{C-\la_\a^2}}$$ and let $x_\be=1+\la_\be^2$,
then (\ref{con}) is equivalent to 
\begin{equation}
 \f{x_\a-1}{x_\a+C-1}+\sum_{\be\neq \a}\Big[\f{C}{C+1-x_\a}-\f{2C^3}{(C+1-x_\a)^2}\f{1}{x_\be
 -\f{x_\a+2C^2-C-1}{x_\a-C-1}}\Big]\leq -1.
\end{equation}

Denote
$$\aligned
\psi(x_\a)=\f{x_\a-1}{x_\a+C-1},&\qquad \varphi(x_\a)=\f{x_\a+2C^2-C-1}{x_\a-C-1},\\
\zeta(x_\a)=\f{C}{C+1-x_\a},&\qquad
\xi(x_\a)=\f{2C^3}{(C+1-x_\a)^2}.
\endaligned
$$

Let
\begin{equation}\label{Om}
\aligned\Om=\big\{&(x_1,\cdots,x_m)\in \R^m: x_\a>C+1,1\leq x_\be<\varphi(x_\a)\text{ for all }\be\neq \a,\g, \\
                    &\qquad x_\g>\varphi(x_\a), \prod_{\be}x_\be=v^2\big\}
                    \endaligned
\end{equation}
and define $f:\Om\ra \R$
$$(x_1,\cdots,x_m)\mapsto \psi(x_\a)+\sum_{\be\neq \a}\Big[\zeta(x_\a)-\f{\xi(x_\a)}{x_\be-\varphi(x_\a)}\Big].$$
We point out that in (\ref{Om}), $\a$ and $\g$ are fixed indices.

Now we claim
\begin{equation}\label{claim}
\sup_\Om f=\sup_\G f
\end{equation}
where
\begin{equation}\aligned
\G=\big\{&(x_1,\cdots,x_m)\in \R^m: x_\a\geq C+1, x_\be=1\text{ for all }\be\neq \a,\g,\\
             &\qquad x_\g\geq \varphi(x_\a),\prod_\be x_\be=v^2\big\}\subset \Om.
\endaligned
\end{equation}
When $m=2$, obviously $\G=\Om$ and (\ref{claim}) is trivial. We
put
$$\varphi_\ep(x_\a)=\varphi(x_\a+\ep),\ \zeta_\ep(x_\a)=\zeta(x_\a+\ep),\ \xi_\ep(x_\a)=\xi(x_\a+\ep)$$
for arbitrary $\ep>0$. If $m\geq 3$, as in the proof of Lemma
\ref{l1}, we define
$$f_\ep=\psi(x_\a)+\sum_{\be\neq \a}\Big[\zeta_\ep(x_\a)-\f{\xi_\ep(x_\a)}{x_\be-\varphi_\ep(x_\a)}\Big],$$
then $f_\ep$ is well-defined on
$$\aligned\Om_\ep=\big\{&(x_1,\cdots,x_m)\in \R^m:x_\a\geq C+1,1\leq x_\be\leq \varphi_{2\ep}(x_\a)\text{ for all }\be\neq \a,\g,\\
                        &\qquad x_\g\geq \varphi_{\f{\ep}{2}}(x_\a),
\prod_{\be}x_\be=v^2\big\}.\endaligned$$
The compactness of $\Om_\ep$ enables us to find $(y_1,\cdots,y_m)\in \Om_\ep$, such that
\begin{equation}\label{mini}
f_\ep(y_1,\cdots,y_m)=\sup_{\Om_\ep} f_\ep.
\end{equation}
Denote $b=\varphi_\ep(y_\a)$, then (\ref{mini}) implies for arbitrary
$\be\neq \a,\g$ that 
$$\f{1}{x_\be-b}+\f{1}{x_\g-b}\geq \f{1}{y_\be-b}+\f{1}{y_\g-b}$$
holds whenever $x_\be x_\g=y_\be y_\g$, $1\leq x_\be\leq
\varphi_{2\ep}(y_\a)$ and $x_\g\geq \varphi_{\f{\ep}{2}}(y_\a)$.
Differentiating both sides yields
$\f{dx_\be}{x_\be}+\f{dx_\g}{x_\g}=0$, thus
\begin{equation}\label{diff}\aligned
d\Big(\f{1}{x_\be-b}+\f{1}{x_\g-b}\Big)&=-\f{dx_\be}{(x_\be-b)^2}-\f{dx_\g}{(x_\g-b)^2}\\
                                       &=\f{(b^2-x_\be x_\g)(x_\g-x_\be)}{(x_\be-b)^2(x_\g-b)^2}\f{dx_\be}{x_\be}.
\endaligned
\end{equation}
Similarly to (\ref{ineq7}), one can prove $y_\a
b^2=\f{y_\a(y_\a+\ep+2C^2-C-1)^2}{(y_\a+\ep-C-1)^2}>9$ when
$\ep_0\leq \ep_1$ (note that $C=1-\ep_0$) and $\ep_1$ is
sufficiently small. In conjunction with $y_\a x_\be x_\g=y_\a
y_\be y_\g\leq v^2<9$, we have $b^2-x_\be x_\g>0$. Hence
(\ref{diff}) implies $y_\be=1$ for all $\be\neq \a, \g$. In other
words, if we put 
$$\aligned
\G_\ep=\big\{(x_1,\cdots,x_m)&\in \R^m: x_\a\geq C+1,\; x_\be=1\text{ for all }\be\neq \a,\g,\\
             &x_\g\geq \varphi_{\f{\ep}{2}}(x_\a),\prod_\be x_\be=v^2\big\},
\endaligned$$
then $\max_{\Om_\ep}f_\ep=\max_{\G_\ep}f_\ep$. Therefore, (\ref{claim}) follows from
 $\Om\subset \bigcup_{\ep>0}\Om_\ep$, $\G\subset \bigcup_{\ep>0}\G_\ep$ and $\lim_{\ep\ra 0}f_\ep=f$.

To prove (\ref{con2}), i.e. $f\leq -1$ , it is sufficient to show
on $\G$,
\begin{equation}
\psi(x_\a)+\zeta(x_\a)-\f{\xi(x_\a)}{\f{v^2}{x_\a}-\varphi(x_\a)}\leq
-1
\end{equation}
whenever $x_\a > C+1$ and $\f{v^2}{x_\a}>\varphi(x_\a)$. After a
straightforward calculation, the above inequality is equivalent to
\begin{equation}\label{con6}
x_\a^3+(2C^2-C-2)x_\a^2+(C^3-3C^2+C+1)x_\a-v^2(x_\a^2-(C+2)x_\a-(C^2-C-1))\geq 0.
\end{equation}
It is easily seen that if
\begin{equation}\label{con5}
\inf_{t^2-(C+2)t-(C^2-C-1)>0}\f{t^3+(2C^2-C-2)t^2+(C^3-3C^2+C+1)t}{t^2-(C+2)t-(C^2-C-1)}> 9.
\end{equation}
then (\ref{con6}) naturally holds and furthermore one can deduce that $IV_\a-\ep_0\big(h_{\a,\a\a}^2+\sum_{\be\neq \a}(h_{\a,\be\be}^2+2h_{\be,\be\a}^2)\big)$ is nonnegative definite.

When $C=1$, (\ref{con5}) becomes
\begin{equation}\label{con7}\inf_{t>\f{3+\sqrt{5}}{2}}\f{t^2(t-1)}{t^2-3t+1}>9.
\end{equation}
If this is true, one can find a positive constant $\ep_3$ to ensure (\ref{con5}) holds true whenever
$\ep_0\leq \ep_3$. Finally,, by taking $\ep_0=\min\{\ep_1,\ep_2,\ep_3\}$ we obtain the final conclusion.

(\ref{con7}) is equivalent to the property that 
$h(t)=t^2(t-1)-9(t^2-3t+1)=t^3-10t^2+27t-9$ has no zeros on $\big(\f{3+\sqrt{5}}{2},+\infty\big)$. $h'(t)=3t^2-20t+27$ implies
$h'(t)<0$ on $\big(\f{3+\sqrt{5}}{2},\f{10+\sqrt{19}}{3}\big)$ and $h'(t)>0$ on $\big(\f{10+\sqrt{19}}{3},+\infty\big)$, hence
$$\inf_{t> \f{3+\sqrt{5}}{2}}h=h\big(\f{10+\sqrt{19}}{3}\big)=\f{187-38\sqrt{19}}{27}>0$$
and (\ref{con7})  follows.

\end{proof}
\bigskip

In conjunction with (\ref{es1}), (\ref{es2}), Lemma \ref{l1} and \ref{l2}, we can arrive at

\begin{thm}\label{thm1}
Let $M^n$ be a submanifold in $\R^{n+m}$ with parallel mean curvature, then for arbitrary $p\in M$ and $P_0\in \grs{n}{m}$,
once $v(\g(p),P_0)\leq 3$, then $\De\big(v(\cdot,P_0)\circ \g\big)\geq 0$ at $p$. Moreover,
if $v(\g(p),P_0)\le
q \be_0<3$, then there exists a positive constant $K_0$, depending only on $\be_0$, such that
\begin{equation}\label{Dv1}
\De\big(v(\cdot,P_0)\circ \g\big)\geq K_0|B|^2
\end{equation}
at $p$.
\end{thm}

We also express this result by saying that the function $v$ satisfying (\ref{Dv1}) is strongly
subharmonic under the condition $v(\g(p),P_0)\leq \be_0<3$.

\begin{rem}
If\, $\log v$ is a strongly subharmonic function, then $v$ is
certainly strongly subharmonic, but the converse is not
necessarily true. Therefore, the above result does not seem to follow 
from Theorem 1.2 in \cite{wang}.
\end{rem}

\Section{Curvature estimates}{Curvature estimates}

Let $z^\a=f^\a(x^1,\cdots,x^n),\a=1,\cdots,m$ be smooth functions
defined on $D_{R_0}\subset \R^n$. Their graph $M=(x,f(x))$ is a
submanifold with parallel mean curvature in $\R^{n+m}$. Suppose
there is $\be_0\in [1,3)$, such that
\begin{equation}\label{slope}
\De_f=\Big[\det\Big(\de_{ij}+\sum_\a \f{\p f^\a}{\p x^i}\f{\p f^\a}{\p x^j}\Big)\Big]^{\f{1}{2}}\leq \be_0.
\end{equation}

Denote by $\eps_1,\cdots,\eps_{n+m}$ the canonical basis of
$\R^{n+m}$ and put $P_0=\eps_1\w\cdots\w\eps_n$. Then by 
(\ref{slope})  $$v(\cdot,P_0)\circ \g\leq \be_0$$ holds
everywhere on $M$. Putting $v=v(\cdot,P_0)\circ \g$,  Theorem
\ref{thm1} tells us
\begin{equation}\label{sub}
\De v\geq K_0(\be_0)|B|^2.
\end{equation}
 Let $\eta$ be a
nonnegative smooth function on $M$ with compact support. Multiplying both sides of (\ref{sub}) by $\eta$ and integrating on $M$
gives
\begin{equation}\label{weak}
K_0\int_M |B|^2 \eta*1\leq -\int_M \n\eta\cdot\n v*1.
\end{equation}

$F:D_{R_0}\mapsto M$ defined by 
$$x=(x^1,\cdots,x^n)\mapsto (x,f(x))$$
is obviously a diffeomorphism. $F_*\f{\p}{\p x^i}=\eps_i+\f{\p f^\a}{\p x^i}\eps_{n+\a}$ implies
$$\big\lan F_*\f{\p}{\p x^i},F_*\f{\p}{\p x^j}\big\ran=\de_{ij}+\sum_\a \f{\p f^\a}{\p x^i}\f{\p f^\a}{\p x^j}.$$
Hence
\begin{equation}
F^* g=\Big(\de_{ij}+\sum_\a \f{\p f^\a}{\p x^i}\f{\p f^\a}{\p x^j}\Big)dx^idx^j
\end{equation}
where $g$ is the metric tensor on $M$. In other words, $M$ is
isometric to the Euclidean ball $D_{R_0}$ equipped with the metric
$g_{ij}dx^i dx^j$  ($g_{ij}=\de_{ij}+\sum_\a \f{\p f^\a}{\p
x^i}\f{\p f^\a}{\p x^j}$). It is easily seen that for arbitrary
$\xi\in \R^n$,
\begin{equation}\label{eig1}
\xi^i g_{ij}\xi^j=|\xi|^2+\sum_\a \Big(\sum_i \f{\p f^\a}{\p x^i}\xi^i\Big)^2\geq |\xi|^2.
\end{equation}
On the other hand, $\De_f\leq \be_0$ implies $\prod_{i=1}^n
\mu_i\leq \be_0^2$, with $\mu_1,\cdots,\mu_n$ the eigenvalues of
$(g_{ij})$, thus
\begin{equation}\label{eig2}
\xi^i g_{ij}\xi^j\leq \be_0^2|\xi|^2\leq 9|\xi|^2.
\end{equation}
In local coordinates, the Laplace-Beltrami operator is 
$$\De=\f{1}{\sqrt{G}}\f{\p}{\p x^i}\Big(\sqrt{G}g^{ij}\f{\p }{\p x^j}\Big).$$
Here $(g^{ij})$ is the inverse matrix of $(g_{ij})$, and $G=\det(g_{ij})=\De_f^2$. From
(\ref{slope}), (\ref{eig1}) and (\ref{eig2}) it is easily seen that
\begin{equation}\label{uniform}
\f{1}{3}|\xi|^2\leq \be_0^{-1}|\xi|^2\leq \xi^i(\sqrt{G}g^{ij})\xi^j\leq \be_0|\xi|^2\leq 3|\xi|^2.
\end{equation}

Following \cite{j} and \cite{j-x-y} we shall make use of the
following abbreviations: For arbitrary $R\in (0,R_0)$, let
\begin{equation}
B_R=\big\{(x,f(x)):x\in D_R\big\}\subset M.
\end{equation}
And for arbitrary $h\in L^\infty(B_R)$ denote
\begin{equation}
\aligned 
&h_{+,R}\mathop{=}\limits^{\text{def.}}\sup_{B_R}h,\qquad 
h_{-,R}\mathop{=}\limits^{\text{def.}}=\inf_{B_R}h,\qquad
\bar{h}_R\mathop{=}\limits^{\text{def.}}\aint{B_R}h=\f{\int_{B_R}h*1}{|\text{Vol}(B_R)|}\\
&|\bar{h}|_{p,R}\mathop{=}\limits^{\text{def.}}\Big(\aint{B_R}|h|^p\Big)^{\f{1}{p}}\; 
(p\in (-\infty,+\infty).
\endaligned
\end{equation}

(\ref{uniform}) shows that $\De$ is a uniform elliptic operator. Moser's Harnack inequality \cite{m} for supersolutions
of elliptic  PDEs in  divergence form gives

\begin{lem}\label{Har}
For a positive superharmonic function $h$ on $B_R$ with $R\in 
(0,R_0]$, $p\in (0,\f{n}{n-2})$ and $\th\in [\f{1}{2},1)$, we 
have the following estimate
$$|\bar{h}|_{p,\th R}\leq \g_1 h_{-,\th R}.$$
Here $\g_1$ is a positive constant only depending on $n$, $p$ and
$\th$, but not  on $h$ and $R$.
\end{lem}

(\ref{sub}) shows the subharmonicity of $v$, and therefore $v_{+,R}-v+\ep$ is a positive superharmonic function on $B_R$
for arbitrary $\ep>0$. With the aid of Lemma \ref{Har}, one can follow \cite{j} to get
\bigskip
\begin{cor}\label{c1}
There is a constant $\de_0\in (0,1)$, depending only on $n$, such that
$$v_{+,\f{R}{2}}\leq (1-\de_0)v_{+,R}+\de_0\bar{v}_{\f{R}{2}}.$$

\end{cor}
\bigskip

Denote by $G^\rho$ the mollified Green function for the  Laplace-Beltrami operator on
$B_R$. Then for arbitrary $p=(y,f(y))\in B_R$, once
$$B_\rho(p)=\big\{(x,f(x))\in M: x\in D_R(y)\big\}\subset B_R$$
we have
$$\int_{B_R}\n G^\rho(\cdot,p)\cdot \n\phi*1=\aint{B_\rho(p)}\phi$$
for every $\phi\in H_0^{1,2}(B_R)$. The apriori estimates for mollified Green functions
of \cite{g-w} tell us

\bigskip
\begin{lem}\label{l4}
With $o:=(0,f(0))$,  we have
\begin{equation}\label{Gr1}\aligned
0\leq G^{\f{R}{2}}(\cdot,o)\leq c_2(n)R^{2-n}&\qquad \text{on }B_R\\
G^{\f{R}{2}}(\cdot,o)\geq c_1(n)R^{2-n}&\qquad \text{on }B_{\f{R}{2}}.
\endaligned
\end{equation}
For arbitrary $p\in B_{\f{R}{4}}$,
\begin{equation}\label{Gr3}
G^\rho(\cdot,p)\leq C(n)R^{2-n}\qquad \text{on }B_R\backslash \bar{B}_{\f{R}{2}}.
\end{equation}
Moreover if $\rho\leq \f{R}{8}$,
\begin{equation}\label{Gr2}
\int_{B_R\backslash \bar{B}_{\f{R}{2}}}\big|\n G^\rho(\cdot,p)\big|^2*1\leq C(n)R^{2-n}.
\end{equation}
\end{lem}
\bigskip

Based on (\ref{weak}), Corollary \ref{c1} and Lemma \ref{l4}, we can
use the method of \cite{j} to derive a  telescoping lemma a la 
Giaquinta-Giusti \cite{g-g} and Giaquita-Hildebrandt \cite{g-h}.
\bigskip

\begin{thm}\label{thm2}
There exists a positive constant $C_1$, only depending on $n$ and $\be_0$, such that for arbitrary $R\leq R_0$,
\begin{equation}\label{tele}
R^{2-n}\int_{B_{\f{R}{2}}}|B|^2*1\leq C_1(v_{+,R}-v_{+,\f{R}{2}})
\end{equation}
Moreover, there exists a positive constnat $C_2$, only depending
on $n$ and $\be_0$, such that for arbitrary $\ep>0$, we can find $R\in
[\exp(-C_2\ep^{-1})R_0,R_0]$, such that
\begin{equation}
R^{2-n}\int_{B_{\f{R}{2}}}|B|^2*1\leq \ep.
\end{equation}

\end{thm}

\begin{proof} With
$$\om^R=R^{-2}\text{Vol}(B_{\f{R}{2}})G^{\f{R}{2}}(\cdot,o)\qquad \text{where }o=(0,f(0)),$$
then
$$\int_{B_R}\n \om^R\cdot \n \phi*1=R^{-2}\int_{B_{\f{R}{2}}}\phi*1.$$
Choosing $(\om^R)^2\in H_0^{1,2}(B_R)$ as a test function in (\ref{weak}), we obtain
$$\aligned
K_0 \int_{B_R}|B|^2(\om^R)^2*1&\leq -\int_{B_R}\n (\om^R)^2\cdot \n v*1=-2\int_{B_R}\n \om^R\cdot \om^R\n(v-v_{+,R})*1\\
                              &=-2\int_{B_R}\n \om^R\cdot\big(\n(\om^R(v-v_{+,R}))-(v-v_{+,R})\n \om^R\big)*1\\
                              &\leq -2\int_{B_R}\n \om^R\cdot \n\big(\om^R(v-v_{+,R})\big)*1\\
                              &=-2R^{-2}\int_{B_{\f{R}{2}}}\om^R(v-v_{+,R})*1.
 \endaligned$$
By (\ref{Gr1}), there exist positive constants $c_3,c_4$, depending only on $n$, such that
$$\aligned
0\leq \om^R\leq c_4 \qquad &\text{on } B_R,\\
\om^R\geq c_3\qquad &\text{on }B_{\f{R}{2}}.
\endaligned$$
Hence
\begin{equation}\label{tele2}\aligned
\int_{B_{\f{R}{2}}}|B|^2*1&\leq -2K_0^{-1}c_4^{-1}c_3^2R^{-2}\int_{B_{\f{R}{2}}}(v-v_{+,R})*1\\
                          &\leq c_5(n,\be_0)R^{n-2}(v_{+,R}-\bar{v}_{\f{R}{2}}).
\endaligned
\end{equation}
By Corollary \ref{c1}, $v_{+,R}-\bar{v}_{\f{R}{2}}\leq \de_0^{-1}(v_{+,R}-v_{+,\f{R}{2}})$; substituting it into
(\ref{tele2}) yields (\ref{tele}).

For arbitrary $k\in \Bbb{Z}^+$, (\ref{tele}) gives
\begin{equation}\aligned
\sum_{i=0}^k (2^{-i}R_0)^{2-n}\int_{B_{2^{-i-1}R_0}}|B|^2*1&\leq C_1\sum_{i=0}^k(v_{+,2^{-i}R_0}-v_{+,2^{-i-1}R_0})\\
                                                          &=C_1(v_{+,R_0}-v_{+,2^{-k-1}R_0})\\
                                                          &\leq C_1(\be_0-1)\leq 2C_1
                                                          \endaligned
\end{equation}
For arbitrary $\ep>0$, we take
$$k=[2C_1\ep^{-1}],$$
then we can find $1\leq j\leq k$, such that
$$(2^{-j}R_0)^{2-n}\int_{B_{2^{-j-1}R_0}}|B|^2*1\leq \f{2}{k+1}C_1\leq \ep.$$
Since $2^{-j}\geq 2^{-k}\geq 2^{-2C_1\ep^{-1}}=\exp\big[-2(\log 2)C_1\ep^{-1}\big]$,
it is sufficient to choose $C_2=-2(\log 2)C_1$.

\end{proof}

\Section{Gauss image shrinking property}{A Gauss image shrinking
property}

\begin{lem}\label{l3}
For arbitrary $a>1$ and $\be_0\in [1,a)$, there exists a positive
constant $\ep_1=\ep_1(a,\be_0)$ with the following property. If
$P_1,Q\in \grs{n}{m}$ satisfies $v(Q,P_1)\leq b\leq \be_0$, then we can find $P_2\in \grs{n}{m}$, such that
$v(P,P_2)\leq a$ for every $P\in \grs{n}{m}$ satisfying $v(P,P_1)\leq b$, and
\begin{equation}
1\leq v(Q,P_2)\leq \left\{\begin{array}{cc}
1 & \text{if }b<\sqrt{2}(1+a^{-1})^{-\f{1}{2}}\\
b-\ep_1 & \text{otherwise.}
\end{array}\right.
\end{equation}
\end{lem}

\begin{proof}
 Obviously $w(P,P)=1$ for every $P\in \grs{n}{m}$, which shows $\grs{n}{m}$ is a submanifold in a Euclidean
 sphere via the Pl\"ucker embedding. Denote by $r(\cdot,\cdot)$ the restriction of the spherical distance on $\grs{n}{m}$,
then by spherical geometry, $w=\cos r$ and hence $v=\sec r$.

Denote $\a=\arccos(a^{-1})$ and $\be=\arccos(b^{-1})$. Now we put $\g=\a-\be$ and
\begin{equation}\label{dis}
 c=\sec \g=(a^{-1}b^{-1}+(1-a^{-2})^{\f{1}{2}}(1-b^{-2})^{\f{1}{2}})^{-1}.
\end{equation}
Once $v(P_2,P_1)\leq c$, the triangle inequality implies
$$r(P,P_2)\leq r(P,P_1)+r(P_2,P_1)\leq \arccos(b^{-1})+\arccos(c^{-1})=\a$$
for every $P$ satisfying $v(P,P_1)\leq b$, and thus $v(P,P_2)\leq a$
 follows.

If $b<\sqrt{2}(1+a^{-1})^{-\f{1}{2}}$, then a direct
calculation shows $\be<\f{\a}{2}$, hence $\g>\be$ and moreover
$v(Q,P_1)\leq b<c$. Thereby $P_2=Q$ is the required point.

Otherwise $b\geq \sqrt{2}(1+a^{-1})^{-\f{1}{2}}$ and hence $c\leq b$. Obviously one of the following two cases
has to occur:

\textit{Case I.} $v(Q,P_1)<c$. One can take $P_2=Q$ to ensure $v(\cdot,P_2)\leq a$ whenever $v(\cdot,P_1)\leq b$. In this case
\begin{equation}\label{below1}
b-v(Q,P_2)=b-1\geq \sqrt{2}(1+a^{-1})^{-\f{1}{2}}-1.
\end{equation}

\textit{Case II.} $v(Q,P_1)\geq c$. Denote by $\th_1,\cdots,\th_m$ the Jordan angles between $Q$ and $P_1$, and put
$L^2=\sum_{1\leq \a\leq m}\th_\a^2$, then as shown in \cite{w}, if we denote the shortest normal geodesic from
$Q$ to $P_1$ by $\g$, then the Jordan angles between $Q$ and $\g(t)$ are $\f{\th_1}{L}t,\cdots,\f{\th_m}{L}t$,
while the Jordan angles between $\g(t)$ and $P_1$ are $\f{\th_1}{L}(L-t),\cdots,\f{\th_m}{L}(L-t)$. Hence
$$\aligned
v(Q,\g(t))&=\prod_\a \sec\big(\f{\th_\a}{L}t\big),\\
v(\g(t),P_1)&=\prod_\a \sec\big(\f{\th_\a}{L}(L-t)\big).
\endaligned$$
Since $t\mapsto \prod_\a \sec\big(\f{\th_\a}{L}(L-t)\big)$ is a strictly decreasing function, there exists
a unique $t_0\in [0,L)$, such that $\prod_\a \sec\big(\f{\th_\a}{L}(L-t_0)\big)=c$. Now we choose
$P_2=\g(t_0)$, then $v(P_2,P_1)=c$ and
\begin{equation}\label{below2}
b-v(Q,P_2)=b-\prod_\a \sec\big(\f{\th_\a}{L}t_0\big).
\end{equation}

It remains to show $b-\prod_\a
\sec\big(\f{\th_\a}{L}t_0\big)$ is bounded from below by a
universal positive constant $\ep_2$. Once this holds true, in conjunction with (\ref{below1}) and (\ref{below2}),
\begin{equation}
\ep_1=\min\{\sqrt{2}(1+a^{-1})^{-\f{1}{2}}-1,\ep_2\}
\end{equation}
is the required constant.

$t_0$ can be regarded as a smooth function on
$$\Om=\big\{(b,\th_1,\cdots,\th_m)\in \R^{m+1},\sqrt{2}(1+a^{-1})^{-\f{1}{2}}\leq b\leq \be_0,
0\leq \th_\a\leq \f{\pi}{2},c\leq \prod_\a \sec(\th_\a)\leq b\big\}$$
which is the unique one satisfying
$$\prod_\a \sec\big(\f{\th_\a}{L}(L-t_0)\big)=c.$$
(By (\ref{dis}), $c$ can be viewed as a function of $b$.) The
smoothness of $t_0$ follows from the implicit function theorem.
Therefore $F:\Om\ra \R$
$$(\th_1,\cdots,\th_m)\mapsto b-\prod_\a \sec\big(\f{\th_\a}{L}t_0\big)$$
is a smooth function on $\Om$. $t_0<L$ implies $F>0$; then the compactness of $\Om$ gives $\inf_\Om F>0$,
and $\ep_2=\inf_\Om F$ is the required constant.

\end{proof}

\begin{rem}
$\ep_1$ is only depending on $a$ and $\be_0$, non-decreasingly
during the iteration process in Theorem \ref{thm4}.
\end{rem}

\begin{thm}\label{thm3}
 Let $M=\big\{(x,f(x)):x\in D_{R_0}\subset \R^n\big\}$ be a graph with parallel mean curvature, and $\De_f\leq \be_0$
 with $\be_0\in [1,3)$. Assume there exists
$P_0\in \grs{n}{m}$, such that $v(\cdot,P_0)\circ \g\leq b$ on $M$ with $1\leq b\leq \be_0$.
If $b<\f{\sqrt{6}}{2}$, then for arbitrary $\ep>0$, one can find a constant $\de\in (0,1)$ depending only on
$n$, $\be_0$ and $\ep$
such that
\begin{equation}\label{es3}
1\leq v(\cdot,P_1)\circ \g\leq 1+\ep\qquad \text{on }B_{\de R_0}
\end{equation}
for a point $P_1\in \grs{n}{m}$.
If $b\geq \f{\sqrt{6}}{2}$, then there are two constants $\de_0\in (0,1)$ and $\ep_1>0$, only depending on $n$ and $\be_0$, such that
\begin{equation}\label{es4}
1\leq v(\cdot,P_1)\circ \g\leq b-\f{\ep_1}{2}\qquad \text{on }B_{\de_0R_0}
\end{equation}
for a point $P_1\in \grs{n}{m}$.
\end{thm}

\begin{proof}

Let $H$ be a smooth function on $\grs{n}{m}$, then $h=H\circ \g$ gives a smooth function on $M$. Let $\eta$ be a
nonnegative smooth function on $M$ with compact support and $\varphi$ be a $H^{1,2}$-function on $M$, then
by Stokes' Theorem,
$$\aligned
0&=\int_M \text{div}(\varphi\eta\n h)*1\\
 &=\int_M \varphi\n\eta\cdot \n h*1+\int_M \eta\n\varphi\cdot \n h*1+\int_M \varphi \eta\De h*1\\
 &=\int_M \varphi\n \eta\cdot \n h*1+\int_M\n\varphi\cdot \n(\eta h)*1-\int_M h\n \varphi\cdot \n\eta*1+\int_M \varphi\eta \De h*1.
\endaligned$$
Hence
\begin{equation}\label{sh0}
 \int_M \n\varphi\cdot \n(\eta h)*1=-\int_M \varphi\n\eta\cdot \n h*1+\int_M h\n\varphi\cdot \n\eta*1-\int_M \varphi\eta\De h*1.
\end{equation}
For arbitrary $R\leq R_0$, we take a cut-off function $\eta$  supported in the interior of
$B_R$, $0\leq \eta\leq 1$, $\eta\equiv 1$ on $B_{\f{R}{2}}$ and $|\n \eta|\leq c_0 R^{-1}$. For every
$\rho\leq \f{R}{8}$, denote by $G^\rho$ the mollified Green function on $B_R$. For arbitrary $p\in B_{\f{R}{4}}$,
inserting $\varphi=G^\rho(\cdot,p)$ into (\ref{sh0}) gives
\begin{equation}\label{sh}\aligned
 &\int_{B_R}\n G^\rho(\cdot,p)\cdot \n(\eta h)*1\\
=&-\int_{B_R}G^\rho(\cdot,p)\n\eta\cdot \n h*1+\int_{B_R}h\n G^\rho(\cdot,p)\cdot \n\eta*1
-\int_{B_R}G^\rho(\cdot,p)\eta\De h*1.\endaligned
\end{equation}
We write (\ref{sh}) as
$$I_\rho=II_\rho+III_\rho+IV_\rho.$$

By the definition of mollified Green functions,
\begin{equation}\label{sh11}
I_\rho=\aint{B_\rho(p)}\eta h=\aint{B_\rho(p)}h.
\end{equation}
Hence
\begin{equation}\label{sh1}
 \lim_{\rho\ra 0^+}I_\rho=h(p).
\end{equation}

Noting that $|d\g|^2=|B|^2$, we have $|\n h|\leq |\n^G H||d\g|=|\n^G H||B|$.
Here and in the sequel, $\n^G$ denotes the Levi-Civita connection on $\grs{n}{m}$. In
conjunction with (\ref{Gr3}), we have
\begin{equation}\label{sh2}\aligned
 |II_\rho|&\leq \int_{T_R}G^\rho(\cdot,p)|\n \eta||\n h|*1\\
          &\leq \sup_{T_R}G^\rho(\cdot,p)\sup_{T_R}|\n \eta|\sup_{\Bbb{V}}|\n^G H|\int_{B_R}|B|*1\\
          &\leq C(n)R^{1-n}\sup_{\Bbb{V}}|\n^G H|\Big(\int_{B_R}|B|^2*1\Big)^{\f{1}{2}}\text{Vol}(B_R)^{\f{1}{2}}\\
          &\leq c_1(n)\sup_{\Bbb{V}}|\n^G H|\Big(R^{2-n}\int_{B_R}|B|^2*1\Big)^{\f{1}{2}}.
\endaligned
\end{equation}
Here $T_R\mathop{=}\limits^{\text{def.}}B_R\backslash \bar{B}_{\f{R}{2}}$ and
$$\Bbb{V}=\{P\in \grs{n}{m}: v(P,P_0)\leq 3\},$$
which is a compact subset of $\Bbb{U}$.

As shown in Section \ref{s1}, there is a one-to-one
correspondence between the points in $\Bbb{U}$ and the 
$n\times m$-matrices. And each $n\times m$-matrix can be viewed
as a corresponding vector in $\R^{nm}$. Define $T:\Bbb{U}\ra
\R^{nm}$
$$Z\mapsto \big(\det(I+ZZ^T)^\f{1}{2}-1\big)\f{Z}{\big(\tr(ZZ^T)\big)^{\f{1}{2}}}$$
Note that
$\big(\tr(ZZ^T)\big)^{\f{1}{2}}=(\sum_{i,\a}Z_{i\a}^2)^{\f{1}{2}}$
equals $|Z|$ when $Z$ is treated as a vector in $\R^{nm}$. Since
$t\in [0,+\infty)\mapsto
\left[\det\big(I+(tZ)(tZ)^T\big)\right]^{\f{1}{2}}$ is a strictly
increasing function and maps $[0,+\infty)$ onto $[1,+\infty)$,
$T$ is a diffeomorphism. By (\ref{v}),  $|T(Z)|=v(P,P_0)-1$. Via
$T$, we can define the mean value of $\g$ on $B_R$ by 
\begin{equation}
\bar\g_R=T^{-1}\Big[\f{\int_{B_R}(T\circ
\g)*1}{\text{Vol}(B_R)}\Big].
\end{equation}
Note that $T$ maps sublevel sets of $v(\cdot,P_0)$ onto Euclidean balls centered at the origin.
Hence the convexity of Euclidean balls gives
\begin{equation}
v(\bar{\g}_R,P_0)\leq \sup_{B_R}v(\cdot,P_0)\circ \g\leq b.
\end{equation}
The compactness of $\Bbb{V}$ ensures the existence of positive constants $K_1$ and $K_2$, such that for arbitrary
$X\in T\Bbb{V}$,
$$K_1|X|\leq |T_* X|\leq K_2|X|.$$

The classical Neumann-Poincar\'{e} inequality says
$$\int_{D_R}|\phi-\bar{\phi}|^2\leq C(n)R^2\int_{D_R}|D \phi|^2.$$
As shown above, $B_R$ can be regarded as $D_R$ equipped with the metric $g=g_{ij}dx^idx^j$, and the eigenvalues
of $(g_{ij})$ are bounded. Hence it is easy to get
$$\int_{B_R}|\phi-\bar{\phi}|^2*1\leq C(n)R^2\int_{B_R}|\n \phi|^2*1.$$
Here $\phi$ can  be a vector-valued function.

Denote by $d_G$ the distance function on $\grs{n}{m}.$ Then, by
using the above Neumann-Poincar\`{e} inequality we have
\begin{equation}\label{poin}\aligned
\int_{B_R}d_G^2(\g,\bar\g_R)*1&\leq K_1^{-2}\int_{B_R}\big|T\circ \g-T(\bar{\g}_R)\big|^2*1\\
                          &\leq C(n)K_1^{-2}R^2\int_{B_R}\big|d(T\circ \g)\big|^2*1\\
                          &\leq C(n)K_1^{-2}K_2^2R^2\int_{B_R}|d\g|^2*1\\
                          &= C(n)K_1^{-2}K_2^2R^2\int_{B_R}|B|^2*1.
                          \endaligned
\end{equation}

Now we write
$$h=H\circ \g=H(\bar{\g}_R)+\big(H\circ \g-H(\bar{\g}_R)\big),$$
then
\begin{equation}\label{sh3}
III_\rho=H(\bar{\g}_R)\int_{B_R}\n G^\rho(\cdot,p)\cdot
\n\eta*1+\int_{T_R}\big(H\circ \g-H(\bar{\g}_R)\big)\n
G^\rho(\cdot,p)\cdot \n \eta*1.
\end{equation}
Similar to (\ref{sh11}),
\begin{equation}\label{sh31}
\lim_{\rho\to 0^+}H(\bar{\g}_R)\int_{B_R}\n G^\rho(\cdot,p)\cdot
\n \eta*1=\lim_{\rho\to
0^+}H(\bar{\g}_R)\aint{B_\rho(p)}\eta=H(\bar\g_R).
\end{equation}
The second term can be controlled by
\begin{equation}\label{sh321}\aligned
&\int_{T_R}\big(H\circ \g-H(\bar{\g}_R)\big)\n G^\rho(\cdot,p)\cdot \n \eta*1\\
\leq&\sup_{\Bbb{V}}|\n^G H|\sup_{T_R}|\n \eta|\int_{T_R}d_G(\g,\bar{\g}_R)|\n G^\rho(\cdot,p)|*1\\
\leq&c_0R^{-1}\sup_{\Bbb{V}}|\n^G H|\Big(\int_{B_R}d_G^2(\g,\bar{\g}_R)\Big)^{\f{1}{2}}
\Big(\int_{T_R}|\n G^\rho(\cdot,p)|^2*1\Big)^{\f{1}{2}}
\endaligned
\end{equation}
Substituting (\ref{Gr2}) and
(\ref{poin}) into (\ref{sh321}) yields
\begin{equation}\label{sh32}
\int_{T_R}\big(H\circ \g-H(\bar{\g}_R)\big)\n
G^\rho(\cdot,p)\cdot \n \eta*1\leq c_2(n)\sup_{\Bbb{V}}|\n^G
H|\Big(R^{2-n}\int_{B_R}|B|^2*1\Big)^{\f{1}{2}}.
\end{equation}

From (\ref{sh}), (\ref{sh1}), (\ref{sh2}), (\ref{sh3}), (\ref{sh31}) and (\ref{sh32}), letting
$\rho\ra 0$ we arrive at
\begin{equation}\label{sh4}
\aligned
h(p)\leq &H(\bar{\g}_R)+c_3(n)\sup_{\Bbb{V}}|\n^G H|\Big(R^{2-n}\int_{B_R}|B|^2*1\Big)^{\f{1}{2}}\\
         &-\limsup_{\rho\ra 0^+}\int_{B_R} G^\rho(\cdot,p)\eta\De h*1.
         \endaligned
\end{equation}
for every $p\in B_{\f{R}{4}}$.

The compactness of $\grs{n}{m}$ implies the existence of a positive constant $K_3$, such that
\begin{equation}\label{sh51}
\big|\n^Gv(\cdot,P)\big|\leq K_3\qquad \text{whenever }1\leq v(\cdot,P)\leq 3
\end{equation}
for arbitrary $P\in \grs{n}{m}$. Hence by
inserting $H=v(\cdot,P)$ into (\ref{sh4}) one can obtain
\begin{equation}\label{sh6}
\aligned
v(\g(p),P)\leq &v(\bar{\g}_R,P)+c_3K_3\Big(R^{2-n}\int_{B_R}|B|^2*1\Big)^{\f{1}{2}}\\
         &-\limsup_{\rho\ra 0^+}\int_{B_R} G^\rho(\cdot,p)\eta\De \big(v(\cdot,P)\circ \g\big)*1.
\endaligned
\end{equation}

By Lemma \ref{l3}, if we put $P_1=\bar{\g}_R$, then $1\leq v(\cdot,P_1)\leq 3$
whenever $1\leq v(\cdot,P_0))\leq b$ provided that $b<\f{\sqrt{6}}{2}$, which implies $v(\cdot,\bar{\g}_R)\circ \g$ is a subharmonic function
on $B_R$. Letting $P=\bar{\g}_R$ in (\ref{sh6}) yields
\begin{equation}\label{sh7}
v(\g(p),\bar{\g}_R)\leq 1+c_3K_3\Big(R^{2-n}\int_{B_R}|B|^2*1\Big)^{\f{1}{2}}
\end{equation}
for all $p\in B_{\f{R}{4}}$. By Theorem \ref{thm2}, for every $\ep>0$, there is $\de\in (0,1)$, depending only
on $n,\be_0$ and $\ep$, such that
\begin{equation}\label{sh8}
R^{2-n}\int_{B_R}|B|^2*1\leq c_3^{-2}K_3^{-2}\ep^2
\end{equation}
for some $R\in [4\de R_0,R_0]$. Substituting (\ref{sh8}) into (\ref{sh7}) gives (\ref{es3}).

If $b\geq \f{\sqrt{6}}{2}$, we put $\ep_1=\ep_1(3,\be_0)$ as given in Lemma \ref{l3}. Then Theorem \ref{thm2} enables us
to find $R\in [4\de_0R_0,R_0]$ such that
\begin{equation}\label{sh53}
R^{2-n}\int_{B_R}|B|^2*1\leq \f{1}{4}c_3^{-2}K_3^{-2}\ep_1^2,
\end{equation}
where $\de_0$ only depends on $n$ and $\be_0$.
Applying Lemma \ref{l3}, one can find $P_1\in \grs{n}{m}$, such that
\begin{equation}\label{sh52}
v(\bar{\g}_R,P_1)\leq b-\ep_1
\end{equation}
and $1\leq v(\cdot,P_1)\leq 3$ whenever $1\leq v(\cdot,P_0)\leq
b$. Theorem \ref{thm1} ensures that $v(\cdot,P_1)\circ \g$ is a
subharmonic function on $B_R$. Taking $P=P_1$
in (\ref{sh6}) yields
\begin{equation}
v\big(\g(p),P_1\big)\leq v(\bar{\g}_R,P_1)+c_3K_3(\f{1}{4}c_3^{-2}K_3^{-1}\ep_1^2)^{\f{1}{2}}\leq
b-\f{\ep_1}{2}
\end{equation}
for all $p\in B_{\f{R}{4}}$. Here we have used (\ref{sh53}) and (\ref{sh52}). From the above inequality (\ref{es4}) immediately follows.

\end{proof}

\Section{Bernstein type results}{Bernstein type results}

Now we can start an iteration as in \cite{h-j-w} and \cite{g-j} to get the following estimates:
\bigskip

\begin{thm}\label{thm4}
Let $M=\big\{(x,f(x)):x\in D_{R_0}\subset \R^n\big\}$ be a graph with parallel mean curvature, and $\De_f\leq \be_0$
 with $\be_0\in [1,3)$, then for arbitrary $\ep>0$, there exists $\de\in (0,1)$, only depending on $n$, $\be_0$ and $\ep$,
 not depending on $f$ and $R_0$, such that
$$1\leq v(\cdot,\g(o))\circ \g\leq 1+\ep\qquad \text{on }B_{\de R_0},$$
where $o=(0,f(0))$. In particular, if $|Df|(0)=0$, then
$$\De_f\leq 1+\ep\qquad \text{on }D_{\de R_0}.$$

\end{thm}

\begin{proof}
Let $\{\eps_1,\cdots,\eps_{n+m}\}$ be canonical orthonormal basis of $\R^{n+m}$ and put $P_0=\eps_1\w\cdots\w\eps_n$.
Then $\De_f\leq \be_0$ implies $v(\cdot,P_0)\leq \be_0$ on $B_{R_0}$. If $\be_0<\f{\sqrt{6}}{2}$,
we put $Q_0=P_0$. Otherwise by Theorem \ref{thm3}, one can find $P_1\in \grs{n}{m}$, such that
\begin{equation}
v(\cdot,P_1)\circ \g \leq \be_0-\ep_1\qquad \text{on }B_{\de_0 R_0}
\end{equation}
with constants $\de_0$ and $\ep_1$ depending only on $n$ and $\be_0$. Similarly for each $j\geq 1$, if
$\be_0-j\ep_1<\f{\sqrt{6}}{2}$, then we put $Q_0=P_j$; otherwise Theorem \ref{thm3} enables us to
 find $P_{j+1}\in \grs{n}{m}$ satisfying
\begin{equation}
v(\cdot,P_{j+1})\circ \g\leq \be_0-(j+1)\ep_1\qquad \text{on }B_{\de_0^{j+1}R_0}.
\end{equation}
Denoting
$$k=\big[(3-\f{\sqrt{6}}{2})\ep_1^{-1}\big]+1,$$
then obviously $\be_0-k\ep_1<\f{\sqrt{6}}{2}$. Hence there exists $Q_0\in \grs{n}{m}$, such that
\begin{equation}
v(\cdot,Q_0)\circ \g\leq b<\f{\sqrt{6}}{2}\qquad \text{on }B_{\de_0^kR_0}.
\end{equation}
Again using Theorem \ref{thm3}, for arbitrary $\ep>0$, there exists $\de_1\in (0,1)$, depending only on $n,\be_0$ and $\ep$, such that
\begin{equation}
v(\cdot,Q_1)\circ \g\leq \sqrt{2}(1+(1+\ep)^{-1})^{-\f{1}{2}}\qquad \text{on }B_{\de_1\de_0^kR_0}
\end{equation}
for a point $Q_1\in \grs{n}{m}$. With $r(\cdot,\cdot)$  as in the proof of Lemma \ref{l3}, then
$$r(\cdot,Q_1)\circ \g=\arccos v(\cdot,Q_1)^{-1}\circ \g\leq \f{1}{2}\arccos (1+\ep)^{-1}.$$
Using the triangle inequality we get
$$r(\cdot,\g(0))\circ \g\leq r(\cdot,Q_1)\circ \g+r(\g(0),Q_1)\circ \g\leq \arccos(1+\ep)^{-1}.$$
Thus $v(\cdot,\g(0))\circ \g\leq 1+\ep$ on $B_{\de_1\de_0^k R_0}$. It is sufficient to put $\de=\de_1\de_0^k$.

\end{proof}

Letting $R_0\ra +\infty$ we can arrive at a Bernstein-type theorem:
\bigskip

\begin{thm}\label{thm5}
Let $z^\a=f^\a(x^1,\cdots,x^n),\ \a=1,\cdots,m$, be smooth functions
defined everywhere in $\R^n$ ($n\geq 3,m\geq 2$). Suppose their graph $M=(x,f(x))$ is a
submanifold with parallel mean curvature in $\R^{n+m}$. Suppose that
there exists a number $\be_0<3$
with
\begin{equation}
\De_f=\Big[\det\Big(\de_{ij}+\sum_\a \f{\p f^\a}{\p x^i}\f{\p f^\a}{\p x^j}\Big)\Big]^{\f{1}{2}}\leq \be_0.\label{be2}
\end{equation}
Then $f^1,\cdots,f^m$ has to be affine linear (representing an affine $n$-plane).
\end{thm}

\noindent{\bf Final remarks} 

\medskip

For any $P_0\in\grs{n}{m}$, denote by $r$ the distance function 
from $P_0$ in $\grs{n}{m}$. The eigenvalues of $\Hess(r)$ were 
computed in \cite{j-x}. Then define
\begin{equation*}
B_{JX}(P_0)=\big\{P\in \grs{n}{m}:\mbox{ sum of any two Jordan
angles between }P\mbox{ and }P_0<\f{\pi}{2}\big\}
\end{equation*} in  the geodesic polar coordinate neighborhood around $P_0$ on the Grassmann
manifold. From (3.2), (3.7) and (3.9)  in \cite{j-x} it turns out 
that $\Hess(r)>0$ on $B_{JX}(P_0).$ Moreover, let $\Si\subset 
B_{JX}(P_0)$ be a closed subset, then $\th_\a+\th_\be\le \be_0
<\f{\pi}2$ and 
$$\Hess(r)\ge  \cot\be_0\ g,$$
where $g$ is the metric tensor on $\grs{n}{m}.$ Hence, the 
composition of the distance function with the Gauss map is a 
strongly subharmonic function on $M$, provided the Gauss image of 
the submanifold $M$ with parallel mean curvature in $\ir{n+m}$ is 
contained in $\Si$. The largest sub-level set of $v(\cdot, P_0)$ 
in $B_{JX}(P_0)$ were studied in \cite{j-x}.  The Theorem 3.2 in 
\cite{j-x} shows that
$$\max\{w(P, P_0);\; P\in\p B_{JX}(P_0)\}=\f{1}{2}.$$ Therefore,
$$\{P\in \grs{n}{m},\; v(\cdot, P_0)<2\}\subset B_{JX}(P_0),$$ and 
$$\{P\in \grs{n}{m};\; v(\cdot, P_0)= 
2\}\bigcap\p B_{JX}(P_0)\neq\emptyset.$$

On the other hand, we can compute directly. From (\ref{He}) we 
also have 
$$\aligned
\Hess(v(\cdot,P_0))&=\sum_{m+1\leq i\leq n,\a}v\ 
\om_{i\a}^2+\sum_{\a}(1+2\la_\a^2)v\ \om_{\a\a}^2
                                          +\sum_{\a\neq \be}\la_\a\la_\be v\ \om_{\a\a}\otimes\om_{\be\be}\\
&\qquad\qquad+\sum_{\a<\be}\Big[(1+\la_\a\la_\be)v\Big(\f{\sqrt{2}}{2}(\om_{\a\be}
+\om_{\be\a})\Big)^2\\
&\hskip2in+(1-\la_\a\la_\be)v\Big(\f{\sqrt{2}}{2}(\om_{\a\be}-\om_{\be\a})\Big)^2\Big].
\endaligned
$$
It follows that $v(\cdot, P_0)$ is strictly convex on 
$B_{JX}(P_0).$ Moreover, if  $\th_\a+\th_\be\le 
\be_0<\f{\pi}{2},$ then
$$\Hess (v(\cdot, P_0))\ge (1-\tan\th_\a\tan\th_\be)v\  g=\f{\cos(\th_\a+\th_\be)}{\cos\th_\a\cos\th_\be}v\,g
\ge\cos\be_0 v\,g$$ where $g$ is the metric tensor of 
$\grs{n}{m}$ and  
$$\De v(\g(\cdot), P_0)\ge \cos\be_0 v|B|^2\ge \cos\be_0 |B|^2.$$
Now, we define
$$\Sigma(P_0)=B_{JX}(P_0)\bigcup\{P\in \grs{n}{m};\;
v(\cdot,P_0)<3\}\subset\grs{n}{m}.$$ The function $v(\cdot, P_0)$ 
is not convex on all of $\Si(P_0)$. But, its precomposition with the 
Gauss map could be a strongly subharmonic  function on $M$ under 
suitable conditions. 

Therefore, we could obtain a more general result: Let $M$ be a 
complete submanifold in $\ir{n+m}$ with parallel mean curvature. 
If its image under the Gauss map is contained in a closed subset 
of $\Sigma(P_0)$ for some $P_0\in \grs{n}{m}$, then $M$ has to be 
an affine linear subspace.

\bigskip\bigskip

\bibliographystyle{amsplain}

\end{document}